\documentclass[12pt]{article}
\usepackage{bbm}
\usepackage{mathrsfs}
\usepackage{amsfonts}
\usepackage{amsmath,amssymb}
\baselineskip 5.5pt \lineskip 5.5pt \numberwithin{equation}{section}

\def\qed{\hfill$\Box$\par}

\def\qed{\ \ \ifhmode\unskip\nobreak\fi\ifmmode\ifinner
         \else\hskip5pt\fi\fi
 \hbox{\hskip5pt\vrule width4pt height6pt depth1.5pt\hskip 1 pt}}
\def\a{\alpha}

\def\cl{\centerline}

\def\vs{\vspace*}

\newtheorem{theo}{Theorem}[section]
\newtheorem{exam}{Example}[section]
\newtheorem{lemm}[theo]{Lemma}

\newtheorem{defi}[theo]{Definition}
\newtheorem{coro}[theo]{Corollary}
\newtheorem{prop}[theo]{Proposition}

\begin{document}
\cl {\large\bf \vs{10pt} Local cocycle $3$-Hom-Lie Bialgebras}
\cl {{\large\bf \vs{10pt} and $3$-Lie Classical Hom-Yang-Baxter Equation }
\noindent\footnote{Supported by the National Science Foundation of
China (No. 11047030) and the Science and Technology
Program of Henan Province (No. 152300410061).\\
* the corresponding author: Yongsheng Cheng, yscheng@ustc.edu.cn}} \vs{6pt}

\cl{ Mengping Wang$^{1)}$, Linli Wu$^{1)}$, Yongsheng Cheng$^{1), 2)*}$}
\cl{$1)$ \small School
of Mathematics and Statistics, Henan
University, Kaifeng 475004, China}
\cl{ $2)$ \small Institute of Contemporary Mathematics, Henan
University, Kaifeng 475004, China}

\vs{16pt}

\noindent{{\bf Abstract.} In this paper, we introduce
$3$-Hom-Lie bialgebras whose compatibility conditions between the multiplication
and comultiplication are given by local cocycle conditions. We study
a twisted $3$-ary version of the Yang-Baxter Equation, called
the $3$-Lie classical Hom-Yang-Baxter Equation ($3$-Lie CHYBE),
which is a general form of $3$-Lie classical Yang-Baxter Equation studied in
\cite{[BGS]} and prove that the bialgebras induced by the solutions of $3$-Lie
CHYBE induce the coboundary local cocycle $3$-Hom-Lie bialgebras.  \vs{5pt}


\vs{6pt}
\par
\cl{\bf\S1. \ Introduction}
\setcounter{section}{1}\setcounter{theo}{0}\setcounter{equation}{0}
\par
$3$-Lie algebra as a generalization of Lie algebra,
has attracted attention from several fields of mathematics and physics.
In particular, the study of $3$-Lie algebras plays an important
role in string theory. For example,
the structure of $3$-Lie algebras is applied to the study of
supersymmetry and gauge symmetry transformations of the
world-volume theory of multiple coincident M2-branes \cite{[BGS], [BCLM], [LCM], [MCL]}.
The classical Yang-Baxter Equation (shorthand for CYBE) has
far-reaching mathematical significance.
The CYBE is closely related to many topics in mathematical physics,
including Hamiltonian structures, Kac-Moody algebras, Poisson-Lie groups,
quantum groups, Hopf algebras, and Lie bialgebras \cite{[D], [Y]}.
There are many known solutions of the CYBE. For example,
there is a 1-1 correspondence between triangular Lie bialgebras and solutions of
the CYBE \cite{[B], [CS2]}.

Motivated by recent work on Hom-Lie bialgebras and the Hom-Yang-Baxter equation \cite{[Y]},
$3$-Lie algebra and $3$-Lie classical Yang-Baxter equation \cite{ [BGS]},
in this paper, we introduce $3$-Hom-Lie bialgebras
whose compatibility conditions between the multiplication
and comultiplication are given by local cocycle condition, and is
called the local cocycle $3$-Hom-Lie bialgebra. We study
a twisted $3$-ary version of the classical Yang-Baxter Equation, called
the $3$-Lie classical Hom-Yang-Baxter Equation (shorthand for $3$-Lie CHYBE),
which is a general form of $3$-Lie classical Yang-Baxter Equation studied in
\cite{[BGS]}, and use the solutions of $3$-Lie
CHYBE to induce the coboundary local cocycle $3$-Hom-Lie bialgebras.

The paper is organized as follow. In Section 2, we define the
$3$-Hom-Lie coalgebra, the $3$-Hom-Lie bialgebra
with the derivation compatibility, the local cocycle $3$-Hom-Lie
bialgebra with cocycle compatibility, and the $3$-Lie CHBYE.
In Section 3, we define the coboundary local cocycle $3$-Hom-Lie bialgebra
and prove that it can be constructed by a multiplicative
$3$-Hom-Lie algebra and a solution of the $3$-Lie CHBYE.

\vs{10pt}
\par

\cl{\bf\S2. \ $3$-Hom-Lie bialgebra with local cocycle condition}
\setcounter{section}{2}\setcounter{theo}{0}\setcounter{equation}{0}
\par
$3$-Hom-Lie algebra as a generalization of $3$-Lie algebra,
its derivations and representation were
introduced in \cite{[LCM]}. For $T=x_{1}\otimes x_{2}\otimes\cdots\otimes x_{n}\in
L^{\otimes n}$ and $1\leqslant i\leqslant j \leqslant n,$
define the $(ij)$-switching operator as
$$ \sigma_{ij}(T)=x_{1}\cdots\otimes x_{j}\otimes\cdots\otimes x_{i}
\otimes\cdots\otimes x_{n}.$$
\begin{defi}
\label{3-hom-lie def}
A $3$-Hom-Lie algebra $(L, [\cdot,\cdot,\cdot],\a)$ is a
vector space $L$ endowed with a $3$-ary linear skew-symmetry multiplication
$\mu: L\otimes L\otimes L \rightarrow L$ satisfying
\begin{eqnarray}
\nonumber\!\!\!\!\!\!\!\!\!\!\!\!&&
\mu(1+\sigma_{12})=0, \mu(1+\sigma_{23})=0,\\
\nonumber\!\!\!\!\!\!\!\!\!\!\!\!&&
\mu(\a\otimes\a\otimes\mu)(1-\omega_{1}-\omega_{2}-\omega_{3})=0,
\end{eqnarray}
where $1$ denotes the identity operator and $\omega_{i}:
L^{\otimes n}\rightarrow L^{\otimes n},$ $1\leq i \leq 3,$
\begin{eqnarray}
\label{3-hom-lie-4}\!\!\!\!\!\!\!\!\!\!\!\!&&
\omega_{1}(x_{1}\otimes x_{2}\otimes x_{3}\otimes x_{4}
\otimes x_{5})=x_{3}\otimes x_{4}\otimes x_{1}\otimes x_{2}\otimes x_{5},\\
\label{3-hom-lie-5}\!\!\!\!\!\!\!\!\!\!\!\!&&
\omega_{2}(x_{1}\otimes x_{2}\otimes x_{3}\otimes x_{4}
\otimes x_{5})=x_{4}\otimes x_{5}\otimes x_{1}\otimes x_{2}\otimes x_{3},\\
\label{3-hom-lie-6}\!\!\!\!\!\!\!\!\!\!\!\!&&
\omega_{3}(x_{1}\otimes x_{2}\otimes x_{3}\otimes x_{4}
\otimes x_{5})=x_{5}\otimes x_{3}\otimes x_{1}\otimes x_{2}\otimes x_{4}.
\end{eqnarray}

If $\alpha[\cdot, \cdot, \cdot]
=[\cdot,\cdot,\cdot]\circ\alpha^{\otimes^{3}},$
we say $L$ is multiplicative.
\end{defi}


For any integer $k$, recall that a linear map
$D: L\rightarrow L$ is called an $\alpha^k$-derivation
of the multiplicative $3$-Hom-Lie algebra $(L, [\cdot, \cdot, \cdot], \alpha)$,
if $D\circ\alpha=\alpha\circ D$ and
$$D[x, y, z]=[D(x), \alpha^k(y), \alpha^k(z)]+[\alpha^k(x), D(y),
\alpha^k(z)]+[\alpha^k(x), \alpha^k(y), D(z)]. $$
For any $x, y\in L$ satisfying $\alpha(x)=x$, $\alpha(y)=y$,
define $ad_k(x, y): L\rightarrow L$ by $ad_k(x, y)(z)=[x, y, \alpha^k(z)]$, $\forall z\in L$.
Then it is easy to prove that $ad_k(x, y)$ is a $\alpha^{k+1}$-derivation
of $(L, [\cdot, \cdot, \cdot], \alpha)$, which we call an inner $\alpha^{k+1}$-derivation.

\begin{defi}
A representation $\rho$ of a multiplicative $3$-Hom-Lie algebra
$(L,[\cdot,\cdot,\cdot],\a)$ on the vector space $V$ with respect to
$A\in gl(V)$, is a linear map $\rho:L\wedge L\rightarrow gl(V)$,
such that, for any $a,b,c,d\in L$,
\begin{eqnarray}
\nonumber\!\!\!\!\!\!\!\!\!\!\!\!&&
(i) \ \rho(\a(a),\a(b))\circ A =A\circ\rho(a,b),\\
\nonumber\!\!\!\!\!\!\!\!\!\!\!\!&&
(ii) \ \rho(\a(b),\a(c))\rho(a,d)+\rho(\a(c), \a(a))\rho(b,d)\\
\nonumber\!\!\!\!\!\!\!\!\!\!\!\!&&-
\rho([a, b, c], \a(d))\circ A+\rho(\a(a),\a(b))\rho(c,d)=0,\\
\nonumber\!\!\!\!\!\!\!\!\!\!\!\!&&
(iii) \ \rho(\a(c),\a(d))\rho(a,b)-\rho(\a(a),\a(b))\rho(c,d)\\
\nonumber\!\!\!\!\!\!\!\!\!\!\!\!&&+
\rho([a,b,c],\a(d))\circ A+\rho(\a(c),[a,b,d])\circ A=0.
\end{eqnarray}
\end{defi}

We often call $(V, A)$ a Hom-module of $(L,[\cdot,\cdot,\cdot],\a)$.

\begin{lemm}
\label{adjo-rep}
Let $(L, [\cdot, \cdot, \cdot], \alpha)$ be a multiplicative
$3$-Hom-Lie algebra.
Define $ad_1: L\wedge L\rightarrow gl(L)$ by
$ad_1(x, y)(z)=[\alpha(x), \alpha(y), z]$, $\forall z\in L$.
Then $ad_1$ is a representation of $(L, [\cdot, \cdot, \cdot], \alpha)$.
\end{lemm}

Give a representation $(\rho, V),$ denote by
$C^{p}_{\a, A}(L, V)$ the set of $p$-cochains
$$C^{p}_{\a,A}(L,V):=\{linear \, maps \, f:\, \underbrace {L\otimes
L\otimes\cdots\otimes L}_{n}\rightarrow V, A\circ f=f\circ \a^{\otimes^{n}}\}.$$
The coboundary operators associated to the Hom-module are given in \cite{[MCL]}.
For $n\geq1,$ the coboundary operator $\delta:C_{\a,A}^{n}(L,V)
\rightarrow C_{\a,A}^{n+2}(L,V)$ is defined as follows:

For $f\in C_{\a,A}^{2n-1}(L,V)$,
\begin{eqnarray}
\nonumber &\!\!\!\!\!\!\!\!\!\!\!\!\!\!\!\!
&\delta_{hom}^{2n-1}f(x_{1}, x_{2}, \cdots,  x_{2n+1})\\[6pt]
\nonumber &\!\!\!\!\!\!\!\!\!\!\!\!\!\!\!\!
&=\rho(\a^{n-1}(x_{2n}), \a^{n-1}(x_{2n+1}))
f(x_{1}, x_{2}, \cdots, x_{2n-1})\\[6pt]
\nonumber &\!\!\!\!\!\!\!\!\!\!\!\!\!\!\!\!
&-\rho(\a^{n-1}(x_{2n-1}), \a^{n-1}(x_{2n+1}))
f(x_{1}, x_{2}, \cdots, x_{2n-2}, x_{2n})\\[6pt]
\nonumber &\!\!\!\!\!\!\!\!\!\!\!\!\!\!\!\!
&+\sum_{k=1}^{n}(-1)^{n+k}\rho(\a^{n-1}(x_{2k-1}), \a^{n-1}(x_{2k}))
f(x_{1}, \cdots, \hat{x_{2k-1}}, \hat{x_{2k}}, \cdots, x_{2n+1})  \\[6pt]
\nonumber &\!\!\!\!\!\!\!\!\!\!\!\!\!\!\!\!
&+\sum_{k=1}^{n}\sum_{j=2k+1}^{2n+1}(-1)^{n+k+1}
f(\a(x_{1}), \cdots, \hat{x_{2k-1}}, \hat{x_{2k}},
\cdots, [x_{2k-1}, x_{2k},x_{j}], \cdots, \a(x_{2n+1})) ;
\end{eqnarray}

for $f\in C_{\a,A}^{2n}(L,V),$
\begin{eqnarray}
\nonumber &\!\!\!\!\!\!\!\!\!\!\!\!\!\!\!\!
&\delta_{hom}^{2n}f(y,x_{1}, x_{2}, \cdots, x_{2n+1})\\[6pt]
\nonumber &\!\!\!\!\!\!\!\! &
=\rho(\a^{n}(x_{2n}), \a^{n}(x_{2n+1}))f(y, x_{1}, x_{2}, \cdots, x_{2n-1})\\[6pt]
\nonumber &\!\!\!\!\!\!\!\!\!\!\!\!\!\!\!\!
&-\rho(\a^{n}(x_{2n-1}), \a^{n}(x_{2n+1}))f(y, x_{1}, x_{2}, \cdots, x_{2n-2}, x_{2n})\\[6pt]
\nonumber &\!\!\!\!\!\!\!\!\!\!\!\!\!\!\!\!
&+\sum_{k=1}^{n}(-1)^{n+k}\rho(\a^{n}(x_{2k-1}),\a^{n}(x_{2k}))
f(y,x_{1}, \cdots, \hat{x_{2k-1}}, \hat{x_{2k}}, \cdots, x_{2n+1}) \\[5pt]
\nonumber &\!\!\!\!\!\!\!\!\!\!\!\!
&+\sum_{k=1}^{n}\sum_{j=2k+1}^{2n+1}(-1)^{n+k+1}f(\a(y),\a(x_{1}), \cdots,
\hat{x_{2k-1}}, \hat{x_{2k}}, \cdots, [x_{2k-1}, x_{2k}, x_{j}], \cdots, \a(x_{2n+1})) .
\end{eqnarray}

Let $L$ be a $3$-Hom-Lie algebra and $(\rho, V)$ be the
representation of $L.$ A linear map $f: \, L\rightarrow V$
is called a $1$-cocycle on $L$
associated to $(\rho, V)$ if it satisfies
\begin{eqnarray}
\nonumber\!\!\!\!\!\!\!\!\!\!\!\!&&
f([x,y,z])=\rho(x,y)f(z)+\rho(y,z)f(x)+\rho(z,x)f(y), \forall x,y,z\in L.
\end{eqnarray}

\begin{defi}\label{def3-hom-lie-co-1}
A $3$-Hom-Lie coalgebra is a triple $(L, \Delta,\a)$ consisting of
a linear space $L$ , a bilinear map $\Delta: L\rightarrow
L\otimes L\otimes L$ and a linear map $\a: L\rightarrow L$
satisfying
\begin{eqnarray}
\nonumber\!\!\!\!\!\!\!\!\!\!\!\!&&
\Delta+\sigma_{12}\Delta=0, \Delta+\sigma_{23}\Delta=0,
\mbox{ \ $($skew-symmetry$)$,}\\
\label{3-hom-lie co-2}\!\!\!\!\!\!\!\!\!\!\!\!&&
(1-\omega_{1}-\omega_{2}-\omega_{3})(\a\otimes\a\otimes
\Delta)\Delta=0.\mbox{ \ $($Hom-coJacobi
identity$)$,}
\end{eqnarray}
where $1, \omega_{1}, \omega_{2}, \omega_{3}:L^{\otimes^{5}}
\rightarrow L^{\otimes^{5}}$ satisfying identities
$(\ref{3-hom-lie-4}),$
$(\ref{3-hom-lie-5}),$ $(\ref{3-hom-lie-6})$ respectively.
\end{defi}

We call $\Delta$ the cobracket. If $\Delta\circ\alpha
=\alpha^{\otimes^{3}}\circ\Delta,$ then we say $L$ is comultiplicative.

Similar with the case of Lie bialgebra, suitable extensions
of these conditions to the context of $3$-Lie algebras are not
equivalent, leading to different extensions of Hom-Lie bialgebra \cite{[BGS]}.
In \cite{[BCLM]},  the authors introduced $3$-Lie bialgebra with
derivation compatibility condition. According to the idea, we give the
definition of $3$-Hom Lie bialgebra with derivation compatibility condition.
\begin{defi}
A $3$-Hom-Lie bialgebra with derivation compatibility condition is a
quadruple $(L,[\cdot,\cdot,\cdot],\Delta,\a)$ such that

(1) $(L,[\cdot,\cdot,\cdot],\a)$ is a $3$-Hom-Lie algebra,

(2) $(L, \Delta, \a)$ is a $3$-Hom-Lie coalgebra, 

(3) $\Delta$ and $[\cdot,\cdot,\cdot]$ satisfy the following
derivation compatibility condition:
\begin{eqnarray}
\label{3-hom-lie bi-1}\!\!\!\!\!\!\!\!\!\!\!\!&&
\Delta[x,y,z]=ad^{(3)}_{1}(x,y)\Delta(z)+ad^{(3)}_{1}(y,z)\Delta(x)+
ad^{(3)}_{1}(z,x)\Delta(y),
\end{eqnarray}
where $ad^{(3)}_{1}(x, y), ad^{(3)}_{1}(y, z), ad^{(3)}_{1}(z, x):
L\otimes L\otimes L\rightarrow L\otimes L\otimes L$ are
$3$-ary linear mappings satisfying, for any $u, v, w\in L,$
\begin{eqnarray}
\nonumber
ad^{(3)}_{1}(x, y)(u\otimes v\otimes w)&\!\!\!\!\!\!\!\! &
=(ad_{1}(x, y)\otimes\alpha\otimes\alpha)(u\otimes v\otimes w)
+(\alpha\otimes ad_{1}(x, y)\otimes\alpha)(u\otimes v\otimes w)\\[6pt]
\nonumber&\!\!\!\!\!\!\!\!\!\!\!\!\!\!\!\!
&+(\alpha\otimes\alpha\otimes ad_{1}(x, y))(u\otimes v\otimes w).
\end{eqnarray}
\end{defi}

Unfortunately, unlike the case of Hom-Lie algebra, it is difficult to
develop the relations between $3$-Hom-Lie bialgebra with derivation
compatibility condition and Classical Hom-Yang-Equation. Similar with \cite{[BGS]},
next, we will introduce $3$-Hom-Lie bialgebra with cocycle
compatibility condition related with a natural
extension of the classical Yang-Baxter equation.

\begin{defi}
A local cocycle $3$-Hom-Lie bialgebra is a
quadruple $(L,[\cdot,\cdot,\cdot],\Delta,\a)$ satisfies:
(1) $(L,[\cdot,\cdot,\cdot],\a)$ is a $3$-Hom-Lie algebra,\\
(2) $(L,\Delta,\a)$ is a $3$-Hom-Lie coalgebra,\\
(3) $\Delta=\Delta_{1}+\Delta_{2}+\Delta_{3}:L\rightarrow
L\otimes L\otimes L$ is a linear map and
the following condition are satisfied:\\
(i) $\Delta_{1}$ is a 1-cocycle associated to
the representation $(L\otimes L\otimes L,ad_{1}\otimes\a\otimes\a),$\\
(ii) $\Delta_{2}$ is a 1-cocycle associated to
the representation $(L\otimes L\otimes L,\a\otimes ad_{1}\otimes\a),$\\
(iii) $\Delta_{3}$ is a 1-cocycle associated
to the representation $(L\otimes L\otimes L,\a\otimes\a\otimes ad_{1}).$\\
\end{defi}

For any $r=\sum_{i}x_{i}\otimes
y_{i}\in L\otimes L,$ define $r_{pq}$ puts $x_{i}$ at the $p$-th position,
$y_{i}$ at the $q$-th position and $1$ elsewhere in an $n$-tensor.
For example, when $n=4,$ we have
$$ r_{12}=\sum_{i}x_{i}\otimes y_{i}\otimes1\otimes1\in
L^{\otimes^{4}}, r_{21}=\sum_{i}y_{i}\otimes x_{i}\otimes1\otimes1\in L^{\otimes^{4}}.$$
Define $[[r,r,r]]^{\a}\in L\otimes^{4}$ by
\begin{eqnarray}
\nonumber
[[r,r,r]]^{\a}&\!\!\!\!\!\!\!\!
&= [r_{12},r_{13},r_{14}]+[r_{12},r_{23},r_{24}]+[r_{13},r_{23},r_{34}]
+[r_{14},r_{24},r_{34}]\\[6pt]
\nonumber &\!\!\!\!\!\!\!\!\!\!\!\!\!\!\!\!
&=\sum_{i,j,k}([x_{i},x_{j},x_{k}]\otimes\a(y_{i})\otimes\a(y_{j})
\otimes\a(y_{k})+\a(x_{i})\otimes[y_{j},x_{j},x_{k}]\otimes\a(y_{j})
\otimes\a(y_{k}) \\[6pt]
\label{3-lie CHYBE 1} &\!\!\!\!\!\!\!\!\!\!\!\!\!\!\!\!
&+\a(x_{i})\otimes\a(x_{j})\otimes[y_{i},y_{j},x_{k}]
\otimes\a(y_{k})+\a(x_{i})\otimes\a(x_{j})\otimes\a(x_{k})
\otimes[y_{i},y_{j},y_{k}]).
\end{eqnarray}

\begin{defi}
\label{3-Lie CHYBE}
Let $L$ be a $3$-Hom-Lie algebra and $r\in L\otimes L.$ The equation
$$[[r,r,r]]^{\a}=0$$
is called the $3$-Lie classical Hom-Yang-Baxter equation $($ $3$-Lie CHYBE $)$.
\end{defi}

The following is the classical Yang-Baxter equation with respect to Lie algebra
\begin{eqnarray}
\nonumber\!\!\!\!\!\!\!\!
&&[[r,r]]^{\a}=[r_{12}, r_{13}]+[r_{12}, r_{13}]+[r_{12}, r_{13}]\\
\nonumber\!\!\!\!\!\!\!\! &
=\!\!\!\!&\sum\limits_{i,j}[x_i,x_j]\otimes\a(y_i)\otimes\a(y_j)+
\sum\limits_{i,k}\a(x_i)\otimes [y_i, x_k]\otimes\a(y_k)+
\sum\limits_{j,k}\a(x_j)\otimes\a(x_k)\otimes[y_j,y_k]\\
\nonumber\!\!\!\!\!\!\!\! &
=\!\!\!\!&0.
\end{eqnarray}
From this, we can see that (\ref{3-lie CHYBE 1}) can be regarded as a natural extension
of the classical Yang-Baxter equation with respect to Lie algebra.

For any $r=\sum_{i}x_{i}\otimes y_{i}\in L\otimes L,$ $x\in L$, define
\begin{eqnarray}
\label{3-lie CHYBE 2}
\left\{\begin{array}{ll}
\Delta_{1}(x)=\sum_{i,j}[x,x_{i},x_{j}]\otimes\a(y_{j})\otimes\a(y_{i}),\\
\Delta_{2}(x)=\sum_{i,j}\a(y_{i})\otimes[x,x_{i},x_{j}]\otimes\a(y_{j}),\\
\Delta_{3}(x)=\sum_{i,j}\a(y_{j})\otimes\a(y_{j})\otimes[x,x_{i},x_{j}].
\end{array}
\right.
\end{eqnarray}
\begin{prop}
\label{lemma 1}
With the above notations and $\a^{\bigotimes^{2}}(r)=r$, we have\\
(i) $\Delta_{1}$ is a $1$-cocycle associated to the representation
$(L\otimes L\otimes L,ad_{1}\otimes\a\otimes\a), $\\
(ii) $\Delta_{2}$ is a $1$-cocycle associated to the representation
$(L\otimes L\otimes L, \a\otimes ad_{1}\otimes\a), $\\
(iii) $\Delta_{3}$ is a $1$-cocycle associated to the representation
$(L\otimes L\otimes L,\a\otimes\a\otimes ad_{1}). $
\end{prop}
\noindent{\it Proof.~}~For all $x,y,z\in L,$ we have
\begin{eqnarray}
\nonumber
\Delta_{1}([x,y,z])&\!\!\!\!\!\!\!\!
&= \sum_{i,j}[[x,y,z],x_{i},x_{j}]\otimes\a(y_{j})\otimes\a(y_{i})\\[6pt]
\nonumber &\!\!\!\!\!\!\!\!\!\!\!\!\!\!\!\!
&=\sum_{i,j}[[x,y,z],\a(x_{i}),\a(x_{j})]\otimes\a^{2}(y_{j})
\otimes\a^{2}(y_{i})\\[6pt]
\nonumber&\!\!\!\!\!\!\!\!\!\!\!\!\!\!\!\!
&=\sum_{i,j}\Big([[x,x_{i},x_{j}],\a(y),\a(z)]+[[y,x_{i},x_{j}],\a(z),\a(x)]
\\[6pt]
\nonumber &\!\!\!\!\!\!\!\!\!\!\!\!\!\!\!\!
&\quad+[[z,x_{i},x_{j}],\a(x),\a(y)]\Big)\otimes\a^{2}(y_{j})
\otimes\a^{2}(y_{i})\\[6pt]
\nonumber&\!\!\!\!\!\!\!\!\!\!\!\!\!\!\!\!
&=(ad_{1}(y,z)\otimes\a\otimes\a)\Delta_{1}(x)+
(ad_{1}(z,x)\otimes\a\otimes\a)\Delta_{1}(y)\\[6pt]
\nonumber &\!\!\!\!\!\!\!\!\!\!\!\!\!\!\!\!
&\quad+(ad_{1}(x,y)\otimes\a\otimes\a)\Delta_{1}(z).
\end{eqnarray}
Therefore, $\Delta_{1}$ is a $1$-cocycle associated to
the representation $(L\otimes L\otimes L,
ad_{1}\otimes\a\otimes\a).$
The other two statements can be proved similarly.
\hfill$\Box$
\vskip16pt

\cl{\bf 3. Coboundary local cocycle $3$-Hom-Lie bialgebra and the $3$-Lie CHBYE}
\setcounter{section}{3}\setcounter{theo}{0}\setcounter{equation}{0}

\begin{defi}
A (multiplicative) coboundary local cocycle $3$-Hom-Lie bialgebra
$(L, [\cdot, \cdot, \cdot], \Delta, \a, r)$ is
a (multiplicative) local cocycle $3$-Hom-Lie bialgebra, and there exists
an element $r=\sum_{i}x_{i}\otimes y_{i}\in L\otimes L$ such that
$\a^{\otimes^{2}}(r)=r$ and $\Delta=\Delta_{1}+\Delta_{2}+\Delta_{3},$
where $\Delta_{1},\Delta_{2},\Delta_{3}$ are induced by $r$
as in (\ref{3-lie CHYBE 2}).
\end{defi}

For $a\in L$ and $1\leq i\leq5,$
define the linear map $\otimes_{i}L\rightarrow\otimes^{5}L$ by
inserting $a$ at the $i$-th position. For example,
for any $t=t_{1}\otimes t_{2}\otimes t_{3}\otimes t_{4},$
we have $t\otimes_{2}a=t_{1}\otimes a\otimes t_{2}\otimes t_{3}\otimes t_{4}.$
The following theorem shows that how we can construct a
multiplicative coboundary local cocycle $3$-Hom-Lie bialgebra.
\begin{theo}
Let $(L,[\cdot,\cdot,\cdot],\a)$ be a multiplicative $3$-Hom-Lie
algebra and $r\in L^{\otimes^{2}}$ satisfying $\a^{\otimes^{2}}(r)=r$ and
\begin{eqnarray}
\nonumber&\!\!\!\!\!\!\!\!
&\sum_{i}(ad_{1}(x_{i},x)\otimes\a\otimes\a\otimes\a\otimes\a)
([[r,r,r]]_{1}^{\a}\otimes_{2}\a(y_{i}))\\[6pt]
\nonumber&\!\!\!\!\!\!\!\!
&+\sum_{i}(\a\otimes ad_{1}(x,x_{i})\otimes\a\otimes\a\otimes\a)
([[r,r,r]]_{1}^{\a}\otimes_{1}\a(y_{i}))\\[6pt]
\nonumber&\!\!\!\!\!\!\!\!
&+\sum_{i}(\a\otimes\a\otimes ad_{1}(x,x_{i})\otimes\a\otimes\a)
([[r,r,r]]_{2}^{\a}\otimes_{5}\a(y_{i}))\\[6pt]
\nonumber&\!\!\!\!\!\!\!\!
&+\sum_{i}(\a\otimes\a\otimes ad_{1}(x_{i},x)\otimes\a\otimes\a)
([[r,r,r]]_{2}^{\a}\otimes_{4}\a(y_{i}))\\[6pt]
\nonumber&\!\!\!\!\!\!\!\!
&+\sum_{i}(\a\otimes\a\otimes\a\otimes ad_{1}(x,x_{i})\otimes\a)
([[r,r,r]]_{2}^{\a}\otimes_{3}\a(y_{i}))\\[6pt]
\nonumber&\!\!\!\!\!\!\!\!
&+\sum_{i}(\a\otimes\a\otimes\a\otimes ad_{1}(x_{i},x)\otimes\a)
([[r,r,r]]_{3}^{\a}\otimes_{5}\a(y_{i}))\\[6pt]
\nonumber&\!\!\!\!\!\!\!\!
&+\sum_{i}(\a\otimes\a\otimes\a\otimes\a\otimes ad_{1}(x,x_{i}))
([[r,r,r]]_{3}^{\a}\otimes_{4}\a(y_{i}))\\[6pt]
\label{***}&\!\!\!\!\!\!\!\!
&+\sum_{i}(\a\otimes\a\otimes\a\otimes\a\otimes ad_{1}(x_{i},x))
([[r,r,r]]_{3}^{\a}\otimes_{3}\a(y_{i}))=0,
\end{eqnarray}
where
$$[[r, r, r]]_{1}^{\a}:=[r_{12}, r_{13}, r_{14}]+[r_{12}, r_{23}, r_{24}]
-[r_{13}, r_{32}, r_{34}]-[r_{14}, r_{42}, r_{43}],$$
$$[[r, r, r]]_{2}^{\a}:=[r_{12}, r_{31}, r_{14}]-[r_{21}, r_{32}, r_{24}]
-[r_{31},r_{32},r_{34}]-[r_{41},r_{42},r_{34}],$$
$$[[r,r,r]]_{3}^{\a}:=-[r_{12}, r_{13}, r_{41}]+[r_{21}, r_{23}, r_{42}]
-[r_{31}, r_{32}, r_{43}]-[r_{41}, r_{42}, r_{43}].$$
Let $\Delta: \, L\rightarrow L^{\otimes^{2}}$ as $\Delta=\Delta_{1}
+\Delta_{2}+\Delta_{3}$ as in (\ref{3-lie CHYBE 2}).
Then $(L, [\cdot, \cdot, \cdot], \Delta, \a, r)$ is a multiplicative
coboundary local cocycle $3$-Hom-Lie bialgebra.
\end{theo}
\noindent{\it Proof.~}~Let $r=\sum_{i}x_{i}\otimes y_{i}.$
First we will prove that  $\Delta=\Delta_{1}+\Delta_{2}
+\Delta_{3}$ commutes with $\a$.
For
$x\in L$, using Definition \ref{3-lie CHYBE 2},
$\a[\cdot,\cdot,\cdot]=[\cdot,\cdot,\cdot]\circ\a^{\otimes^{3}}$
and the assumption $\a^{\otimes^{2}}(r)=r,$ we have
\begin{eqnarray}
\nonumber \!\!\!\!\!\!\!\!\!\!\!\!\Delta_{1}(\a(x))&=&\sum_{i,j}[\a(x),x_{i},x{j}]
\otimes\a(y_{j})\otimes\a(y_{i})\\
\nonumber \!\!\!\!\!\!\!\!\!\!\!\!&
=&\sum_{i,j}[\a(x),\a(x_{i}),\a(x{j})]
\otimes\a^{2}(y_{j})\otimes\a^{2}(y_{i})\\
\nonumber \!\!\!\!\!\!\!\!\!\!\!\!&
=&\a^{\otimes^{3}}\Delta_{1}(x).
\end{eqnarray}
Similarly, $\Delta_{2}(\a(x))=\a^{\otimes^{3}}\Delta_{2}(x)$ and $
\Delta_{3}(\a(x))=\a^{\otimes^{3}}\Delta_{3}(x),$ then we can obtain
$$\Delta(\a(x))=\Delta_{1}(\a(x))+\Delta_{2}(\a(x))
+\Delta_{3}(\a(x))=\a^{\otimes^{3}}\Delta(x). $$

Second, we show that $\Delta$ is anti-symmetric. Since
$$\sigma_{12}\Delta_{1}(x)=\sum_{i,j}\a(y_{j})\otimes[x,x_{i},x_{j}]\otimes \a(y_{i})=\sum_{i,j}\a(y_{i})\otimes[x,x_{j},x_{i}]\otimes \a(y_{j})=-\Delta_{2}(x),$$
$$\sigma_{12}\Delta_{2}(x)=\sum_{i,j}[x,x_{i},x_{j}]\otimes\a(y_{i})
\otimes\a(y_{j})=-\Delta_{1}(x),$$
$$\sigma_{12}\Delta_{3}(x)=\sum_{i,j}\a(y_{i})
\otimes\a(y_{j})\otimes[x,x_{i},x_{j}]=-\Delta_{3}(x).$$
Hence $\sigma_{12}\Delta(x)=-\Delta(x).$ Similarly,
we have $\sigma_{23}\Delta(x)=-\Delta(x).$

Finally, we show that the $3$-Hom-co-Jacobi identity
of $\Delta$ is equivalent to (\ref{***}).
Since each $\Delta$ contains three terms, there are
$36$ terms in $3$-Hom-co-Jacobi identity. Let $G_{i},1\leq i\leq5,$
denote the sum of these terms where $x$ is at the $i$-th
position in the 5-tensors. Thus
$$ G_{1}+G_{2}+G_{3}+G_{4}+G_{5}=0. $$
There are $6$ terms in $G_{1}:$
$$G_{1}=G_{11}+G_{12}+G_{13}+G_{14}+G_{15}+G_{16},$$
where
$$G_{11}=\sum_{i,j,k,l}[[x,x_{i},x_{j}],x_{k},x_{l}]
\otimes\a(y_{l})\otimes\a(y_{k})\otimes\a^{2}(y_{j})
\otimes\a^{2}(y_{i}),$$
$$G_{12}=\sum_{i,j,k,l}[[x,x_{i},x_{j}],x_{k},x_{l}]
\otimes\a(y_{l})\otimes\a^{2}(y_{i})\otimes\a(y_{k})
\otimes\a^{2}(y_{j}),$$
$$G_{13}=\sum_{i,j,k,l}[[x,x_{i},x_{j}],x_{k},x_{l}]
\otimes\a(y_{l})\otimes\a^{2}(y_{j})\otimes\a^{2}(y_{i})
\otimes\a(y_{k}),$$
$$G_{14}=-\sum_{i,j,k,l}\a([x,x_{i},x_{j}])
\otimes\a^{2}(y_{j})\otimes[\a(y_{i}),x_{k},x_{l}]\otimes\a(y_{l})
\otimes\a(y_{k}),$$
$$G_{15}=-\sum_{i,j,k,l}\a([x,x_{i},x_{j}])
\otimes\a^{2}(y_{j})\otimes\a(y_{k})\otimes[\a(y_{i}),x_{k},x_{l}]
\otimes\a(y_{l}),$$
$$G_{16}=-\sum_{i,j,k,l}\a([x,x_{i},x_{j}])
\otimes\a^{2}(y_{j})\otimes\a(y_{l})\otimes\a(y_{l})
\otimes[\a(y_{i}),x_{k},x_{l}].$$
By Hom-Jacobi identity, we have
\begin{eqnarray}
\nonumber
&&G_{11}+G_{12}+G_{13}\\[6pt]
\nonumber\!\!\!\!\!\!\!\!&=& \sum_{i,j,k,l}[[x_{i},x_{j},x_{k}],\a(x),\a(x_{l})]
\otimes\a^{2}(y_{l})\otimes\a^{2}(y_{k})\otimes\a^{2}(y_{j})
\otimes\a^{2}(y_{i})\\[6pt]
\nonumber \!\!\!\!\!\!\!\!\!\!\!\!\!\!\!\!
&=&\sum_{i,j,k,l}(ad_{1}(x_{l},x)\otimes\a\otimes\a\otimes\a\otimes\a)
[x_{k},x_{j},x_{i}]\otimes\a(y_{l})\otimes\a(y_{k})\otimes\a(y_{j})
\otimes\a(y_{i})\\[6pt]
\nonumber\!\!\!\!\!\!\!\!\!\!\!\!\!\!\!\!
&=&\sum_{i}(ad_{1}(x_{l},x)\otimes\a\otimes\a\otimes\a\otimes\a)
[r_{12},r_{13},r_{14}]\otimes_{2}\a(y_{l}).
\end{eqnarray}
Furthermore, we have
\begin{eqnarray}
\nonumber
G_{14}&\!\!\!\!\!\!\!\!
&=\sum_{i,j,k,l}(ad_{1}(x_{j},x)\otimes\a\otimes\a\otimes\a\otimes\a)
\a(x_{i})\otimes\a(y_{j})\otimes[y_{i},x_{l},x_{k}]\otimes\a(y_{l})
\otimes\a(y_{k})\\[6pt]
\nonumber&\!\!\!\!\!\!\!\!\!\!\!\!\!\!\!\!
&=\sum_{j}(ad_{1}(x_{j},x)\otimes\a\otimes\a\otimes\a\otimes\a)
[r_{12},r_{23},r_{24}]\otimes_{2}\a(y_{j}).
\end{eqnarray}
 Similarly,
$$G_{15}=-\sum_{j}(ad_{1}(x_{j},x)\otimes\a\otimes\a\otimes\a\otimes\a)
[r_{13},r_{32},r_{34}]\otimes_{2}\a(y_{j}),$$
$$G_{16}=\sum_{j}(ad_{1}(x_{j},x)\otimes\a\otimes\a\otimes\a\otimes\a)
[r_{14},r_{42},r_{43}]\otimes_{2}\a(y_{j}).$$
Therefore, we obtain
$$G_{1}=\sum_{i}(ad_{1}(x_{i},x)\otimes\a\otimes\a\otimes\a\otimes\a)
[[r,r,r]]_{1}^{\a}\otimes_{2}\a(y_{i}).$$
In a similar manner, we have
$$G_{2}=\sum_{i}(\a\otimes ad_{1}(x_{i},x)\otimes\a\otimes\a\otimes\a)
[[r,r,r]]_{1}^{\a}\otimes_{1}\a(y_{i}).$$

There are $8$ terms in $G_{3}:$
$$G_{3}=G_{31}+G_{32}+G_{33}+G_{34}+G_{35}+G_{36}+G_{37}+G_{38},$$
where
$$G_{31}=\sum_{i,j,k,l}\a(y_{l})\otimes\a(y_{k})\otimes
[[x,x_{i},x_{j}],x_{k},x_{l}]\otimes \a^{2}(y_{j})\otimes \a^{2}(y_{i}),$$
$$G_{32}=-\sum_{i,j,k,l}\a^{2}(y_{j})\otimes \a^{2}(y_{i})\otimes
[[x,x_{i},x_{j}],x_{k},x_{l}]\otimes \a(y_{l})\otimes\a(y_{k}),$$
$$G_{33}=\sum_{i,j,k,l}[\a(y_{j}),x_{k},x_{l}]\otimes \a(y_{l})\otimes
\a[x,x_{i},x_{j}]\otimes \a(y_{k})\otimes \a^{2}(y_{i}),$$
$$G_{34}=\sum_{i,j,k,l}\a(y_{k})\otimes [\a(y_{j}),x_{k},x_{l}]\otimes
\a[x,x_{i},x_{j}]\otimes \a^{2}(y_{j})\otimes \a(y_{k}),$$
$$G_{35}=\sum_{i,j,k,l}\a(y_{l})\otimes \a(y_{k})\otimes
\a[x,x_{i},x_{j}]\otimes [\a(y_{j}),x_{k},x_{l}]\otimes \a^{2}(y_{i}),$$
$$G_{36}=\sum_{i,j,k,l}[\a(y_{i}),x_{k},x_{l}]\otimes \a(y_{l})\otimes
\a[x,x_{i},x_{j}]\otimes \a^{2}(y_{j})\otimes \a(y_{k}),$$
$$G_{37}=\sum_{i,j,k,l}\a(y_{k})\otimes [\a(y_{i}),x_{k},x_{l}]\otimes
\a[x,x_{i},x_{j}]\otimes \a^{2}(y_{j})\otimes \a(y_{l}),$$
$$G_{38}=\sum_{i,j,k,l}\a(y_{l})\otimes \a(y_{k})\otimes
\a[x,x_{i},x_{j}]\otimes \a^{2}(y_{j})\otimes [\a(y_{i}),x_{k},x_{l}].$$
We have
\begin{eqnarray}
\nonumber
&&G_{31}+G_{32}\\[6pt]\!\!\!\!\!\!\!\!
\nonumber&=& \sum_{i,j,k,l}\a^{2}(y_{l})\otimes\a^{2}(y_{k})\otimes
([\a(x),[x_{i},x_{k},x_{l}],\a(x_{j})]\\[6pt]
&&\nonumber+[\a(x),\a(x_{i}),[x_{j},x_{k},x_{l}]])
\otimes\a^{2}(y_{j})\otimes\a^{2}(y_{i})\\[6pt]
\nonumber\!\!\!\!\!\!\!\!\!\!\!\!\!\!\!\!
&=&-\sum_{i,j,k,l}(\a\otimes \a\otimes ad_{1}(x_{j},x)
\otimes\a\otimes\a)\a(y_{l})\otimes\a(y_{k})\otimes
[x_{l},x_{k},x_{i}]\otimes\a(y_{j})\otimes\a(y_{i})\\[6pt]
\nonumber &\!\!\!\!\!\!\!\!\!\!\!\!\!\!\!\!
&-\sum_{i,j,k,l}
(\a\otimes \otimes\a\otimes ad_{1}(x,x_{i})\a\otimes\a)\a(y_{l})\otimes\a(y_{k})\otimes
[x_{l},x_{k},x_{j}]\otimes\a(y_{j})\otimes\a(y_{i})\\[6pt]
\nonumber \!\!\!\!\!\!\!\!\!\!\!\!\!\!\!\!
&=&-\sum_{j}(\a\otimes \a\otimes ad_{1}(x_{j},x)
\otimes\a\otimes\a)[r_{31},r_{32},r_{34}]
\otimes_{4}\a(y_{j})
\\[6pt]
\nonumber &\!\!\!\!\!\!\!\!\!\!\!\!\!\!\!\!
&-\sum_{i}(\a\otimes \a\otimes ad_{1}(x,x_{i})
\otimes\a\otimes\a)[r_{31},r_{32},r_{34}]
\otimes_{5}\a(y_{i}).
\end{eqnarray}
Furthermore, we have
\begin{eqnarray}
\nonumber
G_{33}+G_{36}\!\!\!\!\!\!\!\!&
&= \sum_{i,j,k,l}(\a\otimes \a\otimes ad_{1}(x,x_{i})
\otimes\a\otimes\a)[x_{l},y_{j},x_{k}]\otimes\a(y_{l})
\otimes\a(y_{j})\\[6pt]
\nonumber &\!\!\!\!\!\!\!\!\!\!\!\!\!\!\!\!
&\otimes\a(y_{k})\otimes\a(y_{i})+\sum_{i,j,k,l}(\a\otimes
\a\otimes ad_{1}(x_{j},x)\otimes\a\otimes\a)\a(y_{j})\otimes\a(y_{k})\\[6pt]
\nonumber &\!\!\!\!\!\!\!\!\!\!\!\!\!\!\!\!
&=\sum_{i}(\a\otimes \a\otimes ad_{1}(x,x_{i})\otimes\a\otimes\a)[r_{12},r_{31},r_{14}]
\otimes_{5}\a(y_{i})\\[6pt]
\nonumber &\!\!\!\!\!\!\!\!\!\!\!\!\!\!\!\!
&+\sum_{j}(\a\otimes \a\otimes ad_{1}(x_{j},x)\otimes\a\otimes\a)[r_{12},r_{31},r_{14}]
\otimes_{4}\a(y_{j}),
\end{eqnarray}
and similarly,
\begin{eqnarray}
\nonumber
G_{34}+G_{37}\!\!\!\!\!\!\!\!&
&= -\sum_{i}(\a\otimes \a\otimes ad_{1}(x,x_{i})\otimes\a\otimes\a)[r_{21},r_{32},r_{24}]
\otimes_{5}\a(y_{i})\\[6pt]
\nonumber &\!\!\!\!\!\!\!\!\!\!\!\!\!\!\!\!
&-\sum_{j}(\a\otimes \a\otimes ad_{1}(x_{j},x)\otimes\a\otimes\a)[r_{21},r_{32},r_{24}]
\otimes_{4}\a(y_{j}),
\end{eqnarray}
\begin{eqnarray}
\nonumber
G_{35}+G_{38}\!\!\!\!\!\!\!\!&
&= -\sum_{i}(\a\otimes \a\otimes ad_{1}(x, x_{i})\otimes\a\otimes\a)[r_{41}, r_{42}, r_{34}]
\otimes_{5}\a(y_{i})\\[6pt]
\nonumber &\!\!\!\!\!\!\!\!\!\!\!\!\!\!\!\!
&-\sum_{j}(\a\otimes \a\otimes ad_{1}(x_{j}, x)\otimes\a\otimes\a)[r_{41}, r_{42}, r_{34}]
\otimes_{4}\a(y_{j}).
\end{eqnarray}
Therefore, we obtain
\begin{eqnarray}
\nonumber
G_{3}\!\!\!\!\!\!\!\!&
&= \sum_{i}(\a\otimes \a\otimes ad_{1}(x,x_{i})\otimes\a\otimes\a)[[r,r,r]]_{2}^{\a}
\otimes_{5}\a(y_{i})\\[6pt]
\nonumber &\!\!\!\!\!\!\!\!\!\!\!\!\!\!\!\!
&+\sum_{i}(\a\otimes \a\otimes ad_{1}(x_{i},x)\otimes\a\otimes\a)[[r,r,r]]_{2}^{\a}
\otimes_{4}\a(y_{i}).
\end{eqnarray}
We similarly obtain
\begin{eqnarray}
\nonumber
G_{4}\!\!\!\!\!\!\!\!&
&= \sum_{i}(\a\otimes \a\otimes\a\otimes ad_{1}(x,x_{i})\otimes\a)[[r,r,r]]_{2}^{\a}
\otimes_{3}\a(y_{i})\\[6pt]
\nonumber &\!\!\!\!\!\!\!\!\!\!\!\!\!\!\!\!
&+\sum_{i}(\a\otimes \a\otimes\a\otimes ad_{1}(x_{i},x)\otimes\a)[[r,r,r]]_{3}^{\a}
\otimes_{5}\a(y_{i}).
\end{eqnarray}
\begin{eqnarray}
\nonumber
G_{5}\!\!\!\!\!\!\!\!&
&= \sum_{i}(\a\otimes \a\otimes\a\otimes\a\otimes ad_{1}(x,x_{i}))[[r,r,r]]_{3}^{\a}
\otimes_{4}\a(y_{i})\\[6pt]
\nonumber &\!\!\!\!\!\!\!\!\!\!\!\!\!\!\!\!
&+\sum_{i}(\a\otimes \a\otimes\a\otimes\a\otimes ad_{1}(x_{i},x))[[r,r,r]]_{3}^{\a}
\otimes_{3}\a(y_{i}).
\end{eqnarray}
This completes the proof.\hfill$\Box$

With the notations above, if $r$ is skew-symmetric,
then we can check that
$$[[r,r,r]]_{1}^{\a}=[[r,r,r]]^{\a},
[[r,r,r]]_{2}^{\a}=-[[r,r,r]]^{\a},
[[r,r,r]]_{3}^{\a}=[[r,r,r]]^{\a}.$$
We can obtain an equivalence description of (\ref{***}):
\begin{eqnarray}
\nonumber&\!\!\!\!\!\!\!\!
&\sum_{i}(ad_{1}(x_{i},x)\otimes\a\otimes\a\otimes\a\otimes\a)
([[r,r,r]]^{\a}\otimes_{2}\a(y_{i}))\\[6pt]
\nonumber&\!\!\!\!\!\!\!\!
&+\sum_{i}(\a\otimes ad_{1}(x,x_{i})\otimes\a\otimes\a\otimes\a)
([[r,r,r]]^{\a}\otimes_{1}\a(y_{i}))\\[6pt]
\nonumber&\!\!\!\!\!\!\!\!
&-\sum_{i}(\a\otimes\a\otimes ad_{1}(x,x_{i})\otimes\a\otimes\a)
([[r,r,r]]^{\a}\otimes_{5}\a(y_{i}))\\[6pt]
\nonumber&\!\!\!\!\!\!\!\!
&-\sum_{i}(\a\otimes\a\otimes ad_{1}(x_{i},x)\otimes\a\otimes\a)
([[r,r,r]]^{\a}\otimes_{4}\a(y_{i}))\\[6pt]
\nonumber&\!\!\!\!\!\!\!\!
&-\sum_{i}(\a\otimes\a\otimes\a\otimes ad_{1}(x,x_{i})\otimes\a)
([[r,r,r]]^{\a}\otimes_{3}\a(y_{i}))\\[6pt]
\nonumber&\!\!\!\!\!\!\!\!
&+\sum_{i}(\a\otimes\a\otimes\a\otimes ad_{1}(x_{i},x)\otimes\a)
([[r,r,r]]^{\a}\otimes_{5}\a(y_{i}))\\[6pt]
\nonumber&\!\!\!\!\!\!\!\!
&+\sum_{i}(\a\otimes\a\otimes\a\otimes\a\otimes ad_{1}(x,x_{i}))
([[r,r,r]]^{\a}\otimes_{4}\a(y_{i}))\\[6pt]
\nonumber&\!\!\!\!\!\!\!\!
&+\sum_{i}(\a\otimes\a\otimes\a\otimes\a\otimes ad_{1}(x_{i},x))
([[r,r,r]]^{\a}\otimes_{3}\a(y_{i}))=0.
\end{eqnarray}
Summarizing the above discussion, we obtain
\begin{coro}
\label{coro1}
Let $L$ be a $3$-Hom-Lie algebra, $\a^{\otimes^{2}}(r)=r$
and $r\in L\otimes L$ skew-symmetric. If
$$[[r,r,r]]^{\a}=0,$$
and $\Delta=\Delta_{1}+\Delta_{2}+\Delta_{3}:
L\rightarrow L\otimes L\otimes L,$ in which $\Delta_{1},
\Delta_{2},\Delta_{3}$ are included by $r$ as in Eq.(\ref{3-lie CHYBE 2}).
Then $(L,[\cdot,\cdot,\cdot],\Delta,\a)$ is a local cocycle $3$-Hom-Lie bialgebra.
\end{coro}

Corollary \ref{coro1} can be regarded as a $3$-Hom-Lie algebra
analogue of the face that $\a^{\otimes^{2}}(r)=r$ and a
skew-symmetric solution of the classical Yang-Baxter
equation gives a local cocycle $3$-Hom-Lie bialgebra.

The next result shows that given a  $3$-Lie algebra endomorphism,
each classical $r$-matrix induces an infinite family of solutions of the $3$-Lie CHYBE.

\begin{theo}
Let $L$ be a $3$-Lie algebra, $r\in L^{\otimes^{2}}$ be a
solution of $3$-Lie CYBE, and $\a:L\rightarrow L$ be a
$3$-Lie algebra endomorphism. Then for each integer
$n\geq0,$ $(\a^{\otimes ^{2}})^{n}(r)$ is a solution
of $3$-Lie CHYBE in the  $3$-Hom-Lie algebra $L_{\a}
=(L, [\cdot, \cdot,\cdot]_{\a}=\a\circ[\cdot,\cdot,\cdot],\a).$
\end{theo}
\noindent{\it Proof.~}~ We can proof that $L_{\a}$
is a $3$-Hom-Lie algebra (in fact, the Hom-Jacobi
identity for $[\cdot,\cdot,\cdot]_{\a}$  is $\a^{3}$
applied to the Jacobi identity of $[\cdot,\cdot,\cdot].$ )
It remains to show that $(\a^{\otimes ^{2}})^{n}(r)$
satisfies the  $3$-Lie CHYBE in the  $3$-Hom-Lie algebra $L_{\a},$ i.e.
$$[[(\a^{\otimes ^{2}})^{n}(r),(\a^{\otimes ^{2}})^{n}(r),
(\a^{\otimes ^{2}})^{n}(r)]]^{\a}=0.$$
Write $r=\sum_{i}x_{i}\otimes y_{i}.$ Using
$\a([\cdot,\cdot,\cdot])=[\cdot,\cdot,\cdot]\circ
\a^{\otimes ^{2}}$ and  the definition
$[\cdot,\cdot,\cdot]_{\a}=\a([\cdot,\cdot,\cdot]),$ we have
\begin{eqnarray}
\nonumber&\!\!\!\!\!\!\!\!
&[[(\a^{\otimes ^{2}})^{n}(r),(\a^{\otimes ^{2}})^{n}(r),
(\a^{\otimes ^{2}})^{n}(r)]]^{\a}\\[6pt]
\nonumber&\!\!\!\!\!\!\!\!
&=\sum_{i,j,k}([\a^{n}(x_{i}),\a^{n}(x_{j}),
\a^{n}(x_{k})]_{\a}\otimes \a(\a^{n}(y_{i}))
\otimes \a(\a^{n}(y_{j}))\otimes \a(\a^{n}(y_{k}))\\[6pt]
\nonumber&\!\!\!\!\!\!\!\!
&+\a(\a^{n}(x_{i}))\otimes[\a^{n}(y_{i}),\a^{n}(x_{j}),
\a^{n}(x_{k})]_{\a}\otimes \a(\a^{n}(y_{j}))\otimes \a(\a^{n}(y_{k}))\\[6pt]
\nonumber&\!\!\!\!\!\!\!\!
&+\a(\a^{n}(x_{i}))\otimes\a(\a^{n}(x_{j}))\otimes[\a^{n}(y_{i}),\a^{n}(y_{j}),
\a^{n}(x_{k})]_{\a}\otimes \a(\a^{n}(y_{k}))\\[6pt]
\nonumber&\!\!\!\!\!\!\!\!
&+\a(\a^{n}(x_{i}))\otimes\a(\a^{n}(x_{j}))\otimes\a(\a^{n}(x_{k}))
\otimes[\a^{n}(y_{i}),\a^{n}(y_{j}),
\a^{n}(y_{k})]_{\a})\\[6pt]
\nonumber&\!\!\!\!\!\!\!\!
&=\a^{n+1}([[r,r,r]])\\[6pt]
\nonumber&\!\!\!\!\!\!\!\!
&=0.
\end{eqnarray}\hfill$\Box$
\begin{exam}
Let $L$ be  the (unique) non-trivial $3$-dimensional complex
$3$-Lie algebra whose non-zero product with respect to a basis
${e_{1},e_{2},e_{3}}$ is given by
$[e_{1},e_{2},e_{3}]=e_{1}. $
If $r=\sum_{i< j}^{3}r_{ij}(e_{i}\otimes e_{j}-e_{j}\otimes e_{i})
\in L\otimes L$ is skew-symmetric and  $\a$ is the $3$-Lie algebra morphism given by
$$\a(e_{1})=a_{11}e_{1},\a(e_{2})=a_{12}+a_{22}e_{2}+a_{32}e_{3},
\a(e_{3})=a_{13}+a_{23}e_{2}+a_{33}e_{3},\forall a_{ij}\in \mathbb{C}(1\leq i,j\leq3)$$
with
$$a_{22}a_{33}-a_{23}a_{32}=1,$$
$$(r_{12}a_{11}-r_{23}a_{13})a_{22}+(r_{13}a_{11}+r_{23}a_{12})a_{23}=r_{12},$$
$$(r_{12}a_{11}-r_{23}a_{13})a_{32}+(r_{13}a_{11}+r_{23}a_{12})a_{33}=r_{13}.$$
Then $L_{\a}=(L,[\cdot,\cdot,\cdot]_{\a}
=\a\circ[\cdot,\cdot,\cdot],\a)$ is a $3$-Hom-Lie algebra and $\a^{\otimes^{2}}(r)=r.$
Moreover, $r$ is a solution of the $3$-Lie CHYBE in $L,$
the corresponding local cocycle $3$-Hom-Lie bialgebra is given by
$$\Delta_{i}(e_{1})=(-1)^{i}r_{23}r\otimes_{i}a_{11}e_{1},$$
$$\Delta_{i}(e_{2})=(-1)^{i+1}r_{13}r\otimes_{i}a_{11}e_{1},$$
$$\Delta_{i}(e_{3})=(-1)^{i}r_{12}r\otimes_{i}a_{11}e_{1}$$
and the comultiplication $\Delta:L\rightarrow \wedge ^{3}L$ is given by
$$ \Delta(e_{1})=-r_{23}^{2}a_{11} e_{1}\wedge e_{2}\wedge e_{3},
\Delta(e_{2})=r_{13}r_{23}a_{11} e_{1}\wedge e_{2}\wedge e_{3},
\Delta(e_{3})=-r_{12}r_{23}a_{11} e_{1}\wedge e_{2}\wedge e_{3},$$
where $e_{1}\wedge e_{2}\wedge e_{3}=\sum_{\sigma\in S_{3}}
sgn(\sigma)e_{\sigma(1)}\otimes e_{\sigma(2)}\otimes e_{\sigma(3)}$
and $S_{3}$ is the permutation group on $\{1,2,3\}.$
\end{exam}
\begin{exam}
Let $L$ be  the (unique) non-trivial $4$-dimensional complex $3$-Lie algebra
whose non-zero product with respect to a basis ${e_{1},e_{2},
e_{3},e_{4}}$ is given by $[e_{1},e_{2},e_{3}]=e_{1}$.
If $r=-e_{3}\otimes e_{4}+e_{4}\otimes e_{3}+\sum_{i<j}^{4}(e_{i}
\otimes e_{j}-e_{j}\otimes e_{i})\in L\otimes L$ is skew-symmetric and
$\a$ is the $3$-Lie algebra morphism given by
$$\a(e_{1})=a_{11}e_{1},$$
$$\a(e_{2})=a_{12}e_{1}+a_{22}e_{2}+e_{3}+e_{4},$$
$$\a(e_{3})=a_{13}e_{1}+a_{23}e_{2}+2e_{3}+2e_{4},$$
$$\a(e_{4})=0,$$
where
$$a_{11}(a_{22}+a_{23})+a_{12}a_{23}-a_{22}a_{13}=1,$$
$$3a_{11}+2a_{12}-a_{13}=1,$$
$$2a_{22}-a_{23}=1.$$
Then $L_{\a}=(L,[\cdot,\cdot,\cdot]_{\a}=
\a\circ[\cdot,\cdot,\cdot],\a)$ is a $3$-Hom-Lie algebra
and $\a^{\otimes^{2}}(r)=r$. Moreover, $r$ is a solution of the
$3$-Lie CHYBE in $L$, the corresponding nontrivial
local cocycle $3$-Hom-Lie bialgebra is given by
$$\Delta_{i}(e_{1})=(-1)^{i}r\otimes_{i}
(a_{11}e_{1}),$$
$$\Delta_{i}(e_{2})=(-1)^{i+1}r\otimes_{i}
(a_{11}e_{1}),$$
$$\Delta_{i}(e_{3})=(-1)^{i}r\otimes_{i}
(a_{11}e_{1}),$$
where $i=1, 2, 3$
and the comultiplication $\Delta:L\rightarrow \wedge ^{3}L$ is given by\\
$$\Delta(e_{1})
=\Delta(e_{3})=-\Delta(e_{2}), \Delta(e_{4})=0, $$
$$\Delta(e_{2})=a_{11}e_{1}\wedge e_{2}\wedge e_{3}+
a_{11}e_{1}\wedge e_{2}\wedge e_{4}, $$
where $e_{1}\wedge e_{2}\wedge e_{3}=\sum_{\sigma\in S_{3}}
sgn(\sigma)e_{\sigma(1)}\otimes e_{\sigma(2)}\otimes e_{\sigma(3)}$
and $S_{3}$ is the permutation group on $\{1,2,3\}.$
\end{exam}

\end{document}